

\input amstex
\documentstyle{amsppt}
\nologo
\NoBlackBoxes

\hcorrection{19mm}

\topmatter
\title  The classification of exceptional Dehn surgeries on 2-bridge knots 
\endtitle 
\rightheadtext{Exceptional surgeries on 2-bridge knots}
\author  Mark Brittenham and Ying-Qing Wu
\endauthor
\address Department of Mathematics, University of North Texas,
     Denton, Texas 76203 \endaddress
\email britten\@unt.edu  \endemail
\address Department of Mathematics, University of
     Iowa, Iowa City, IA 52242 \endaddress 
\email  wu\@math.uiowa.edu
\endemail 

\abstract We will classify all exceptional Dehn surgeries on 2-bridge
knots according to whether they produce reducible, toroidal, or
Seifert fibered manifolds.  
\endabstract 
\endtopmatter
 
\document
\define\aaa{\alpha}

\define\bdd{\partial}
\define\proof{\demo{Proof}}
\define\endproof{\enddemo}

\redefine\hat{\widehat}

\define\rrr{\gamma}

\redefine\bdd{\partial}
\define\em{\it}
\define\Int{\text{Int}}

\head 1. Introduction  \endhead

A nontrivial Dehn surgery on a hyperbolic knot $K$ in $S^3$ is {\it
exceptional\/} if the resulting manifold is either reducible, or
toroidal, or a Seifert fibered manifold whose orbifold is a sphere
with at most three exceptional fibers, called a {\it small Seifert
fibered space}.  Thus an exceptional Dehn surgery is non-hyperbolic,
and using a version of Thurston's orbifold theorem proved by Boileau
and Porti [BP], it can be shown that a non-exceptional surgery on a
2-bridge knot is also hyperbolic, see Remark 6.3.  There has been many
researches in determining how many exceptional surgeries a knot can
admit.  See the survey articles [Go, Lu] for details.

The purpose of this paper is to classify all exceptional Dehn
surgeries on 2-bridge knots.  Hatcher and Thurston [HT] have shown
that there is no reducible surgery on hyperbolic 2-bridge knots.  We
will determine all toroidal surgeries and small Seifert fibered
surgeries on these knots, which will then complete the classification.

We use $[b_1, \ldots, b_n]$ to denote the partial fraction
decomposition $1/(b_1 - 1/(b_2 - \ldots - 1/b_n )\ldots )$.  Recall
that a 2-bridge knot $K$ is non-hyperbolic if and only if $K =
K_{1/q}$ for some $q$, in which case $K$ is a $(2,q)$ torus knot, and
surgery on $K$ is well understood.  $K$ is a {\it twist knot\/} if it
is equivalent to some $K_{p/q}$ with $p/q = [b, \pm 2]$ for some
integer $b$.  Since $[b,\pm2] = [b\mp 1, \mp 2]$, we may assume that
$b$ is even.  Let $K(\rrr)$ be the manifold obtained by $\rrr$ surgery
on $K$.  We always assume that $\rrr \neq \infty$, that is, the
surgery is nontrivial.  With respect to the standard
meridian-longitude pair on $\bdd N(K)$, each slope $\rrr$ is
identified with a rational number, see Rolfsen's book [R].  The
following is the main theorem of this paper.  Note that part (4) of
the theorem is the case of surgery on Figure 8 knot, and is due to
Thurston [Th1].  It is included in the theorem for the sake of
completeness.

\proclaim{Theorem 1.1}
Let $K$ be a hyperbolic 2-bridge knot.

(1) If $K \neq K_{[b_1,b_2]}$ for any $b_1,b_2$, then $K$ admits no
exceptional surgery.

(2) If $K = K_{[b_1,b_2]}$ with $|b_1|,|b_2| > 2$, then $K(\rrr)$ is
exceptional for exactly one $\rrr$, which yields a toroidal manifold.
When both $b_1$ and $b_2$ are even, $\rrr = 0$.  If $b_1$ is odd and
$b_2$ is even, $\rrr = 2b_2$.

(3)  If $K = K_{[2n, \pm 2]}$ and $|n| > 1$,  $K(\rrr)$ is exceptional
for exactly five $\rrr$:  $K(\rrr)$ is toroidal for $\rrr = 0, \mp
4$, and is small Seifert fibered for $\rrr = \mp 1, \mp 2, \mp 3$.

(4) If $K= K_{[2,-2]}$ is the Figure 8 knot, $K(\rrr)$ is exceptional
for only nine $\rrr$: $K(\rrr)$ is toroidal for $\rrr = 0, 4, -4$, and
is Seifert fibered for $\rrr =-1,-2,-3,1,2,3$.  \endproclaim

We will use a result of Hatcher and Thurston [HT] to determine all
toroidal surgeries on a hyperbolic 2-bridge knot $K$, see Lemma 2.2
below.  In general it is more difficult to determine small Seifert
fibered surgeries, due to the fact that there is no essential surfaces
in such a manifold, unless its Euler number is 0.  Here we will use
essential lamination theory developed by Gabai and Oertel [GO].  The
readers are referred to [GO] for definitions and basic properties
concerning essential branched surfaces and essential laminations,
which play a central role in the proof of the theorem.  We will use
Brittenham's criteria [Br], which says that if $M$ is a small Seifert
fibered space containing an essential branched surface $\Cal F$, then
each component of $M - \Int N(\Cal F)$ is an $I$-bundle over some
compact surface $G$.  The idea of the proof is to construct essential
laminations in surgered manifolds whose complementary regions are not
$I$-bundles.  We will apply some techniques developed by Delman [De1,
De2] and Roberts [Ro] to construct essential laminations in the
surgered manifolds.

\medskip

{\em Acknowledgment.}  We are grateful to Alan Reid for some
very useful comments to an earlier version of this paper.  We also
thank Steve Bleiler for pointing out that the Orbifold Theorem and
the Snappea computer program, two tools used in an earlier version
of the paper, are not considered fully proved by many people.  The
present version has avoided using them in the proof.

\head 2. Reducible, toroidal, or Seifert fibered surgeries 
\endhead

Reducible surgeries and toroidal surgeries on 2-bridge knots are
completely determined by Lemmas 2.1 and 2.2.  Certain surgeries on
twist knots are shown to be Seifert fibered in Corollary 2.4.

\proclaim{Lemma 2.1}  {\rm (Hatcher-Thurston)}  Let $K$ be a 2-bridge
knot.  Then $K(r)$ is reducible if and only if $K$ is a $(2, q)$
torus knot, and $r = 2q$.
\endproclaim

\proof See [HT, Theorem 2].
\qed
\enddemo

\proclaim{Lemma 2.2}
Let $K$ be a hyperbolic 2-bridge knot.   

(1)  If $K(\rrr)$ is toroidal for some $\rrr$, then $K = K_{[b_1,
b_2]}$ for some $b_1, b_2$.

(2)  If $|b_i| > 2$ for $i = 1,2$, there is exactly one such
$\rrr$.  When both $b_i$ are even, $\rrr = 0$.  When $b_1$ is odd
and $b_2$ is even, $\rrr = 2b_2$.

(3)  If $K = K_{[2n,2]}$ and $|n| > 1$, $K(\rrr)$ is toroidal if and
only if $\rrr = 0$ or $-4$.  For $K = K_{[2n,-2]}$, $\rrr = 0$ or
$4$.  

(4)  If $K = K_{[2,-2]}$, then $K(\rrr)$ is toroidal if and
only if $\rrr = 0$, $4$, or $-4$. \endproclaim

\proof
We refer the reader to [HT] for notations.  If $K(\rrr)$ is
toroidal, there is an essential punctured torus $T$ in the knot
exterior. By Theorem 1 of [HT], $T$ is carried by some
$\Sigma [b_1, \ldots, b_k]$, where $[b_1, \ldots, b_k]$ is an
expansion of $p/q$.  By the proof of Theorem 2 of [HT], we have
$0 = 2 - 2g = n(2-k)$, where $g$ is the genus of $T$, and $n$ is
the intersection number between $\bdd T$ and a meridian of $K$. 
Therefore $k =2$.  This proves (1).  The rest follows by determining
all the possible expansions of type $[b_1, b_2]$ for $p/q$.  The
boundary slopes of the surfaces can be calculated using Proposition 2
of [HT].  By the proof of [Pr, Corollary 2.1], an incompressible
punctured torus $T$ in the exterior of a 2-bridge knot will become an
essential torus after surgery  along the slope of $\bdd T$. \qed
\enddemo

\proclaim{Lemma 2.3}
Let $L=k_1 \cup k_2$ be the Whitehead link, which is the 2-bridge link
associated to the rational number $p'/q' = [2,2,-2]$.  Let $L(\rrr_1,
\rrr_2)$ be the manifold obtained by $\rrr_i$ surgery on $k_i$.  If
$\rrr_1 = -1/n$ and $\rrr_2 = -1 $, $-2$ or $-3$, then $L(\rrr_1,
\rrr_2)$ is a small Seifert fibered space.
\endproclaim

\proof
By definition $L(\infty,\rrr_2)$ is the manifold obtained from $S^3$
by $\rrr_2$ surgery on $k_2$.  After $-1$ surgery on $k_2$, the knot
$k_1$ becomes a trefoil knot in $L(\infty, -1) = S^3$.  Since the
exterior of a torus knot is a Seifert fibered space with orbifold a
disk with two cones, it is easy to see that all surgeries but one
yield Seifert fibered spaces, each having an orbifold a
disk with at most three cone points.  For this trefoil, the
exceptional surgery has coefficient $-6$, yielding a reducible
manifold.  Thus $L(-1/n, -1)$ is a small Seifert fibered space for
any $n$.  

After $-2$ surgery on $k_2$, the knot $k_1$ becomes a knot in $\Bbb R
P^3 = L(\infty,-2)$.  The link $L$ is drawn in Figure 1(a), where the
curve $C$ is a curve on $\bdd N(k_2)$ of slope $-2$, so it bounds a
disk in $L(\infty,-2)$.  Thus a band sum of $k_1$ and $C$ forms a
knot $k'_1$ isotopic to $k_1$ in $L(\infty,-2)$.  The link $L' = k'_1
\cup k_2$ is shown in Figure 1(b).  Using Kirby Calculus one can show
that $L(-1/n, -2) = L'(-2-1/n, -2)$.  The exterior of $k'_1$ in $S^3$
is a Seifert fibered space with orbifold a disk with two cones, in
which $k_2$ is a singular fiber of index 3.  Thus after $-2$ surgery
on $k_2$, the manifold $L(\infty,-2) - \Int N(k'_1)$ is still Seifert
fibered, with orbifold a disk with two cones.  The fiber slope on
$\bdd N(k'_1)$ is 6.  It follows that all but the 6 surgery on $k'_1$
in $L(\infty,-2)$ yield small Seifert fibered manifolds.  In
particular, $L(-1/n,-2) = L'(-2-1/n, -2)$ are small Seifert fibered
manifolds for all $n$.  

The proof for $\rrr_2 = -3$ is similar.  One can show that the band
sum of $k_1$ and the curve $C$ of slope $-3$ on $\bdd N(k_2)$ is
isotopic to the curve $k'_1$ shown in Figure 1(c), which is a
$(3,-2)$ torus knot. By the same argument as above one can show that 
$L(-1/n, -3) = L(-3-1/n, -3)$ are small Seifert fibered manifolds for
all $n$.   \qed \enddemo

\input epsf.tex

\bigskip
\leavevmode

\epsfxsize=5in
\centerline{\epsfbox{Fig6.ai}}
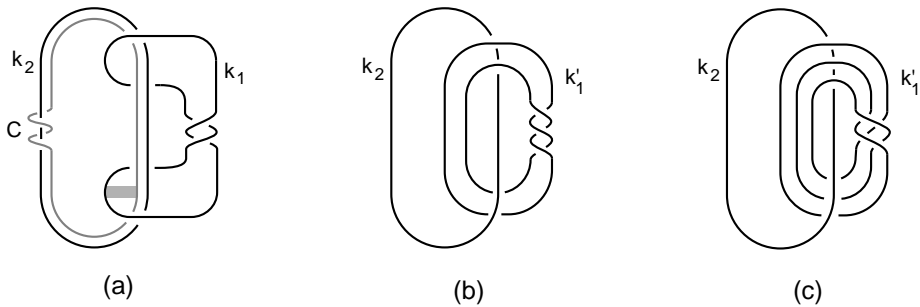
\centerline{Figure 1}

\bigskip

Recall that a 2-bridge knot $K$ is a twist knot if $K = K_{p/q}$,
and $p/q = [2n, \pm 2]$ for some $n$.

\proclaim{Corollary 2.4}
If $K=K_{p/q}$ is a twist knot with $p/q = [2n,\pm 2]$, then 
$K(\rrr)$ is a small Seifert fibered space for
$\rrr = \mp 1$, $\mp 2$ and $\mp 3$.  
\endproclaim

\proof
Consider the case $p/q = [2n,2]$.  The proof for $p/q = [2n,-2]$ is 
similar.  Let $L=k_1 \cup k_2$ be a 2-bridge link associated to
the rational number $p'/q' = [2,2,-2]$.  Notice that after $-1/n$
surgery on $k_1$, the knot $k_2$ becomes the knot $K=K_{[2n,2]}$ in
$S^3 = L(-1/n, \infty)$.  Therefore by Lemma 6, $K(\rrr) =
L(-1/n,\rrr)$ are  small Seifert fibered spaces for $\rrr = -1$, $-2$
and $-3$.  
\qed
\enddemo

\head 3.  Delman's construction, and the proof of Theorem 1.1(1)
\endhead

For each rational number $p/q$, there is associated a diagram
$D(p/q)$, which is the minimal subdiagram of the Hatcher-Thurston
diagram [HT, Figure 4] that contains all minimal paths from $1/0$ to
$p/q$.  See [HT, Figure 5] and [De1].   $D(p/q)$ can be constructed
as follows.  Let $p/q = [a_1, \ldots, a_k]$ be a continued fraction
expansion of $p/q$.  To each $a_i$ is associated a ``fan''
$F_{a_i}$ consisting of $a_i$ simplices, see Figure 2(a) and 2(b) for
the fans $F_4$ and $F_{-4}$.  The edges labeled $e_1$ are called
initial edges, and the ones labeled $e_2$ are called terminal
edges.  The diagram $D(p/q)$ can be constructed by gluing the
$F_{a_i}$ together in such a way that the terminal edge of $F_{a_i}$
is glued to the initial edge of $F_{a_{i+1}}$.  Moreover, if $a_i
a_{i+1} < 0$ then $F_{a_i}$ and $F_{a_{i+1}}$ have one edge in
common, and if $a_i a_{i+1} > 0$ then they have a 2-simplex in
common.  See Figure 2(c) for the diagram of $[2,-2,-4,2]$.  Notice
that the fans $F_{-2}$ and $F_{-4}$ in the figure share a common
triangle.

\bigskip
\leavevmode

\epsfxsize=5in
\centerline{\epsfbox{Fig1.ai}}
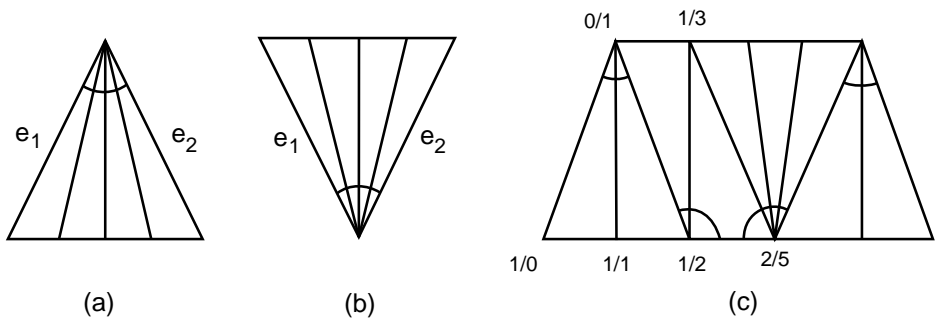
\centerline{Figure 2}

\bigskip

To each vertex $v_i$ of $D(p/q)$ is associated a rational number
$r_i/s_i$.  It has one of the three possible parities:  odd/odd,
odd/even, or even/odd, denoted by $o/o$, $o/e$, and $e/o$,
respectively.  Note that the three vertices of any simplex in
$D(p/q)$ have mutually different parities.

We consider $D(p/q)$ as a graph on a disk $D$, with all vertices on
$\bdd D$, containing $\bdd D$ as a subgraph. The boundary of $D$
forms two paths from the vertex $1/0$ to the vertex $p/q$.  The one
containing the vertex $0/1$ is called the {\it top path}, and the one
containing the vertex $1/1$ is called the {\it bottom path}.  Edges
on the top path are called {\it top edges}.  Similarly for {\it
bottom edges}. 

Let $\Delta_1, \Delta_2$ be two simplices in $D(p/q)$ with an edge in
common.  Assume that the two vertices which are not on the common
edge are of parity $o/o$.  Then the arcs indicated in Figure 3(a)
and (b) are called {\it channels}.  A {\it path\/}  $\aaa$ in
$D(p/q)$ is a union of arcs, each of which is either an edge of
$D(p/q)$ or a channel.

\bigskip
\leavevmode

\epsfxsize=4.5in
\centerline{\epsfbox{Fig2.ai}}
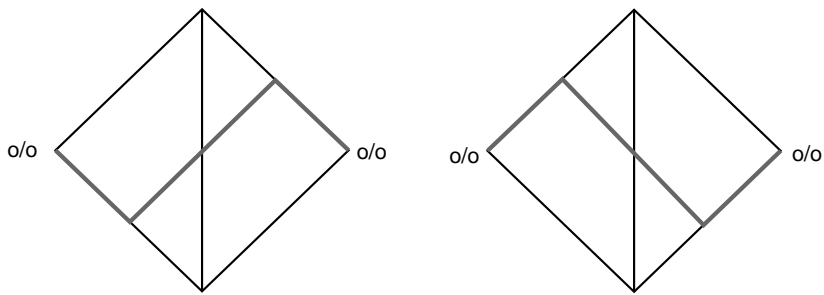
\centerline{Figure 3}

\bigskip

Let $v$ be a vertex on a path $\aaa$ in $D(p/q)$.  Let $e_1, e_2$ be
the edges of $\aaa $ incident to $v$.  Then the {\em corner number\/}
of $v$ in $\aaa$, denoted by $c(v;\aaa)$ or simply $c(v)$, is defined
as the number of simplices in $D(p/q)$ between the edges $e_1$ and
$e_2$.  A path $\aaa$ from $1/0$ to $p/q$ is an {\it allowable path\/}
if it has at least one channel, and $c(v) \geq 2$ for all $v$ in
$\aaa$.

Now assume that $K = K_{p/q}$ is a 2-bridge knot.  Then $q$ is an
odd number.   Recall that $K_{p/q} = K_{p'/q}$ if $p' \equiv p^{\pm
1}$ mod $q$, and $K_{-p/q}$ is the mirror image of $K_{p/q}$.  We
may assume without loss of generality that $p$ is even, and $1<p<q$. 
This is because $K_{(q-p)/q}$ is equivalent to the mirror image of
$K_{p/q}$, so the result of $\gamma$ surgery on the first is the
same as that of $-\gamma$ surgery on the second. Note that $q-p$ and
$p$ have different parity, since $q$ is odd.  The following result is
due to Delman.  See [De1] and [De2, Proposition 3.1].  

\proclaim{Theorem 3.1} {\rm (Delman)} Given an allowable path $\aaa$
of $D(p/q)$, there is an essential branched surface $\Cal F$ in $S^3 -
K$ which remains essential after all nontrivial surgeries on the
knot $K$.  \qed \endproclaim

\proclaim{Lemma 3.2}  If there is an allowable path $\aaa$ in
$D(p/q)$ such that $c(v) > 2$ for some vertex $v$ in $\aaa$, then
$K(\rrr)$ is not a small Seifert fibered space for any $\rrr$.
\endproclaim

\proof
It was shown in [Br, Corollary 4] that if $\Cal F$ is an essential
branched surface in a small Seifert fibered space $M$, then each
component of $M - \Int N(\Cal F)$ is an $I$-bundle over a compact
surface $G$, such that  the vertical surface $\bdd_v N(\Cal F)$
(also called cusps) is the $I$-bundle over $\bdd G$.  It has been
shown in [De1] that for each vertex $v$ of $\aaa$ there is a component
$W_v$ of $S^3-\Int N(\Cal F)$ such that $W_v$ is a solid torus whose
meridian disk intersects the cusps $c(v)$ times. In particular, if
$c(v) > 2$ then $W_v$ is not an $I$ bundle as above.  Since $\Cal F$
is an essential branched surface in $K(\rrr)$, it follows that
$K(\rrr)$ is not a small Seifert fibered space. \qed \enddemo

\proclaim{Lemma 3.3}
Suppose $p$ is even, $q$ is odd, and $1<p<q-1$.
If $p/q$ does not have partial fraction decomposition of type
$[r_1,r_2]$, then $D(p/q)$ has an allowable path $\aaa$ such that
some vertex $v$ on $\aaa$ has $c(v) > 2$.
\endproclaim

\proof
Let $[a_1, \ldots, a_n]$ be the partial fraction decomposition of
$p/q$ such that all $a_i$ are even.  Then $a_1 \geq 2$.  If $a_i = 2$
for all $i$, then $p/q = (q-1)/q$, contradicting our assumption. 
Thus either some $a_i< 0$, or some $a_i > 4$.  We separate the two
cases.

CASE 1.  {\it Some $a_i < 0$.}

Let $a_i$ be the first negative number.  Then $a_{i-1} > 0$, so
there is a sign change.  By [De2] there is a  channel $\aaa_0$ in
$F_{a_{i-1}} \cup F_{a_i}$ starting at a bottom edge and ending at a
top edge, where $F_{a_i}$ is the fan in $D(p/q)$ corresponding to
$a_i$.  Let $\aaa_1$ be the part of the bottom path of $D(p/q)$ from
the vertex $1/0$ to the initial point of $\aaa_0$, and let $\aaa_2$
be the part of the top path from the end point of $\aaa_0$ to the
vertex $p/q$.  Then $\aaa = \aaa_1 \cup \aaa_0 \cup \aaa_2$ is an
allowable path in $D(p/q)$. We need to show that if $c(v) = 2$ for all
vertices $v$ on this path, then $p/q=[r_1,r_2]$ for some $r_1, r_2$.

Consider the vertices on $\aaa_1$.  Since $c(v_i)=2$ for all $v_i$,
each vertex $v_i$ is incident to exactly one non boundary edge
$e_i$ of $D(p/q)$, which must have the other end on a vertex
$v'_i$ in the top path.  If some of these $v'_i$ are different,
then since all faces of $D(p/q)$ are triangles, it is clear that
some $v_j$ on $\aaa_1$ would have at least two non boundary edges,
which would be a contradiction.  Similarly, each vertex on $\aaa_2$
has a unique non boundary edge, leading to a common vertex on the
bottom path, so the diagram $D(p/q)$ looks exactly as in Figure
4(a).  It is the union of two fans $F_{r_1}$ and $F_{r_2}$ with 
$r_1 > 0$, and $r_2 < 0$.    Therefore, $p/q = [r_1,r_2]$.

CASE 2.  {\it Some $a_i \geq 4$.}

 In this case there is a channel $\aaa_0$ with both ends on
the bottom path.  Construct an allowable path $\aaa = \aaa_1 \cup
\aaa_0 \cup \aaa_2$ with $\aaa_1, \aaa_2$ in the bottom path. 
Similar to Case 1, it can be shown that each vertex on $\aaa_i$
has a unique non boundary edge leading to a common vertex $v'_i$ on
the top path, so $D(p/q)$ looks like that in Figure 4(b).  In this
case $p/q = [r_1,r_2]$, with both $r_i >0$.  \qed
\enddemo

\bigskip
\leavevmode

\epsfxsize=4.5in
\centerline{\epsfbox{Fig3.ai}}
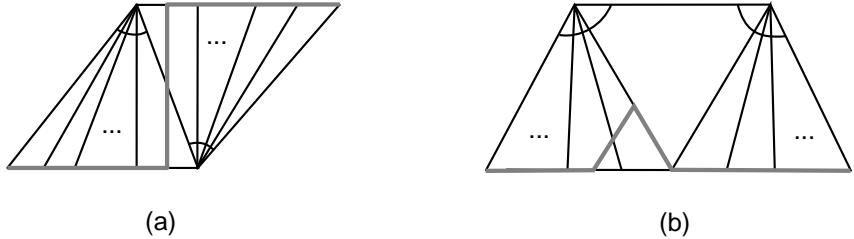
\centerline{Figure 4}

\bigskip

\proclaim{Corollary 3.4}  Let $K$ be a 2-bridge knot.
If $K \neq K_{[b_1,b_2]}$ for any $b_1,b_2$, then $K(\rrr)$ is
non-exceptional for all $\rrr$.
\endproclaim

\proof
$K$ is not a $(2,q)$ torus knot, otherwise $K = K_{[2q\pm 1]}
= K_{[2q, \mp1]}$.  Hence by Lemmas 2.1 and 2.2, $K(\rrr)$ is
irreducible and atoroidal.  By Lemmas 3.3 and 3.2, $K(\rrr)$ is not a
small Seifert fiber space. Therefore, $K(\rrr)$ is non-exceptional.
\qed
\endproof

\head 4. Surgery on twisted Whitehead links,  proof of Theorem
1.1(2) \endhead

A twisted Whitehead link is a two bridge link $L$ associated to a
rational number $[2, r, -2]$ for some $r \neq 0$.  See Figure 5 for
a twisted Whitehead link with $r = -6$.  When $r = \pm 2$, $L$ is a
Whitehead link.  It has been determined exactly which 2-bridge link
complements contain persistent laminations [Wu].  The next lemma
shows that if $|r|>2$ then there is a persistent lamination with
some desired property.

Recall that a slope $\rrr$ of a knot $K$ is an integral slope if it
intersects the meridian of $K$ exactly once.

\proclaim{Lemma 4.1}  Let $L = k_1 \cup k_2$ be a twisted Whitehead
link associated to the rational number $p/q = [2, r, -2]$.  Let
$L(\rrr_1, \rrr_2)$ be the manifold obtained by $\rrr_i$ surgery on
$k_i$.  If $|r| > 2$, and  one of the $\rrr_i$ is not an integral
slope, then $L(\rrr_1, \rrr_2)$ is not a small Seifert fiber space.
\endproclaim

\proof
By considering the mirror image of $L$ if necessary we may assume
that $r < 0$.  If $r$ is even, then $[2,r, -2]$ is a partial
fraction decomposition with even coefficient, and $r \leq -4$. 
There is an allowable path in $D(p/q)$ with two channels, as shown
in Figure 6(a), where $r = -4$.  If $r$ is odd, then $p/q =
[2,r+1, 2]$, in which case $D(p/q)$ also has an allowable path
with two channels.  See  Figure 6(b) for the case $r = -3$.  

\bigskip
\leavevmode

\epsfxsize=2.5in
\centerline{\epsfbox{Fig4.ai}}
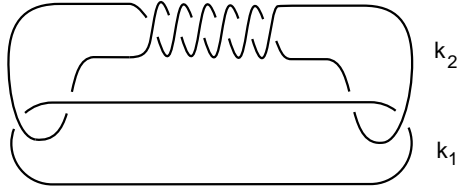
\centerline{Figure 5}

\bigskip

Let $\Cal F$ be the essential branched surface in the link exterior
associated to the above allowable path in $D(p/q)$, as constructed
in [De2].  There is one solid torus component $V_i$ in $S^3 - \Int
N(\Cal F)$ for each $k_i$, containing $k_i$ as a central curve.  
From the construction of $\Cal F$ one can see that each channel
contributes two cusps, one on each $\bdd V_i$.  Actually from [De2,
Figure 3.5] we see that the two cusps corresponding to a channel are
around two points of $L$ on a level sphere with same orientation. 
Since each $k_i$ intersects the sphere at two points with different
orientations, those two cusps must be around different components of
$L$.  One is referred to [Wu] for more details about surgery on
2-bridge links.

\bigskip
\leavevmode

\epsfxsize=4in
\centerline{\epsfbox{Fig5.ai}}
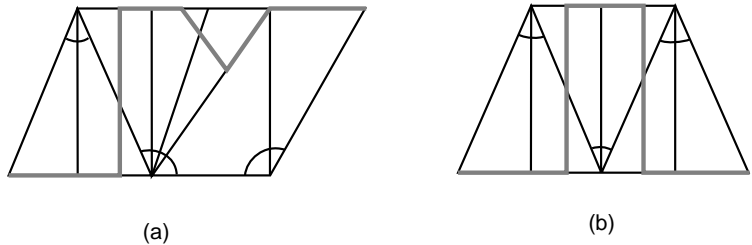
\centerline{Figure 6}

\bigskip

As the allowable path above has two channels, each $V_i$ has two
meridional cusps.  Thus $\Cal F$ remains an essential branched
surface after surgery on $L$.  Moreover, since one of the $\rrr_i$
is non-integral, after surgery $V_i$ becomes a solid torus whose
meridional disk intersects the cusps at least four times.  By [Br,
Corollary 4], the surgered manifold is not a small Seifert fiber
space.  \qed \enddemo

\proclaim{Corollary 4.2} Let $K= K_{[b_1, b_2]}$ be a two bridge knot
with $|b_i| > 2$ for $i=1,2$.  Then $K(\rrr)$ is exceptional for only
one $\rrr$, which yields toroidal manifold.  When both $b_1$ and $b_2$
are even, $\rrr = 0$.  If $b_1$ is odd and $b_2$ is even, $\rrr =
2b_2$.  \endproclaim

\proof
By Lemmas 2.1 and 2.2, $K(\rrr)$ is irreducible, and it is toroidal
for exactly one $\rrr$ as described in the corollary.  So it remains
to show that $K(\rrr)$ is never a small Seifert fiber space.

Since $K$ is a knot, at least one of the $b_i$ is an even number.  We
may assume without loss of generality that $b_1 = 2n$ for some integer
$n$, because $K_{[b_1, b_2]}$ is equivalent to $K_{[b_2,b_1]}$, by
turning the standard diagram for the first knot upside down.  

Let $L = k_1\cup k_2$ be a 2-bridge link associated to the rational
number $p/q = [2,b_2,-2]$.  Notice that after $-1/n$ surgery on $k_1$,
the other component $k_2$ becomes the knot $K = K_{[2n,b_2]}$. 
Therefore, doing $\rrr$ surgery on $K$ is the same as doing $-1/n$
surgery on $k_1$, then doing some $\rrr'$ surgery on $k_2$.  Since 
$-1/n$ is non integral, and $|b_2|>2$, the result follows from Lemma
4.1. \qed
\enddemo

\proclaim{Corollary 4.3}  Let $K =  K_{[b, \pm
2]}$ with $|b| > 2$.  If $\rrr$ is a non integral slope, then
$K(\rrr)$ is not a small Seifert fiber space.
\endproclaim

\proof
As above, $K(\rrr) = L(\rrr, \pm 1)$, where $L = L_{[2, b, -2]}$. 
Since $\rrr$ is non integral, the result follows from Lemma 4.1.
The result also follows from [Br].
\qed
\enddemo

\head 5. Roberts' construction of essential branched surfaces
\endhead

In [Ro] Roberts constructed branched surfaces in certain knot
complements, which can be extended to essential branched surfaces in
$K(\rrr)$ for all $\rrr$ in an infinite interval.  We will describe
her results and construction in this section, and apply them to
surgery on twist knots in the next section.

Let $E(K) = S^3 - \Int N(K)$ be the exterior of a knot $K$ in
$S^3$, let $R'$ be a (possibly non orientable) compact surface in
$S^3$ with $\partial R' = K$.  Let $R = R' \cap E(K)$.

\bigskip
\leavevmode

\epsfxsize=4in
\centerline{\epsfbox{Fig7.ai}}
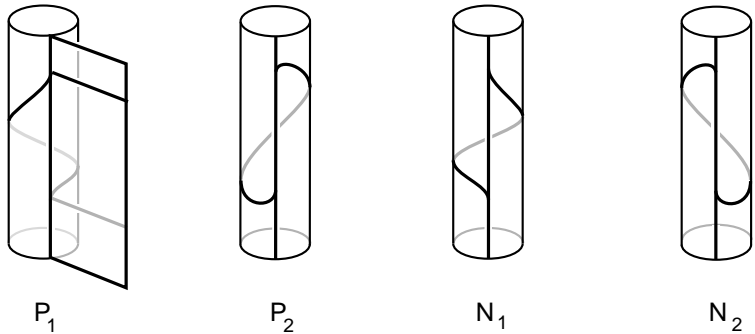
\centerline{Figure 7}

\bigskip

Let $S$ be a surface in $E(K)$ which has interior disjoint from
$R$, and has a single boundary curve $\partial S = a_1 \cup b_1 \cup
\ldots \cup a_n \cup b_n$, where $b_i$ are mutually disjoint arcs on
$R$, and $a_i$ are arcs on $T = \partial N(K)$.  By specifying a
cusp at each $b_i$, the union of $R$ and $S$ becomes a branched
surface $B = \left< R, S\right>$ in $E(K)$.  The cusps will be
assigned in  such a way that each $a_i$ on $T$ is one of the four
types indicated in Figure 7.  Note that $\bdd B$ is a train track on
$T$, and each component of $T - \bdd B$ is a digon, i.e a disk with
two cusps.

\remark{Remark}  The pictures on Figure 7 are mirror images of that
in [Ro, Figure 22].  Thus for example, type $P_1$ here is of type
$N_1$ in [Ro].  Apparently we are using different coordinate
systems.  This paper adapts the convention that the
meridian-longitude pair $(m, l)$ on $T$ is chosen so that when $K$
is endowed with the same orientation as that of $l$, the linking
number $lk(m, K) = 1$, measured using the right hand rule.  See [R]. 
With this convention, types $P_1, P_2$ in Figure 7 will have
positive contributions to any slope $\rrr$ carried in the train
track.  
\endremark

Let $N(B)$ be a regular neighborhood of the branched surface $B$ in
$E(K)$ with the natural $I$-bundle structure.  A {\it surface of
contact\/} is a properly embedded compact surface $P$ in
$N(B)$, transverse to the $I$-fibers, with $\bdd P \subset
\bdd_vN(B) \cup T$, such that the intersection of $\bdd P$ with
each component of the vertical surfaces  $\bdd _v N(B)$ is either
empty or a single arc.

Let $p_1, p_2, n_1, n_2$ be the numbers of $a_i$ of type $P_1, P_2,
N_1, N_2$ respectively.  Let $r$ be the slope of $\bdd R$ on $T$. 
Let 
$$ J = \{ r + (p_1 - n_1) \frac {x}{x+1} + (p_2 - n_2) x \, \,
\vert \, \, x > 0 \}$$
Then Roberts' theorem [Ro, Theorem 2.3] can be stated as 

\proclaim{Theorem 5.1} {\rm (Roberts) } 
If $B = \left< R, S \right>$ constructed above is an essential
branched surface in $E(K)$, and has no planar surface of contact,
then $B$ extends to an essential branched surface $B_{\rrr}$ in
$K(\rrr)$ for all slope $\rrr \in J$.
\endproclaim

The construction of the extended branched surface is as follows. 
Let $T\times I$ be a small neighborhood of $T$ in $E(K)$ with $T =
T\times 0$, such that $B \cap (T\times I) = \bdd B \times I$.  Add the
digons $T \times 1 - B$ to $B$, and branched so that the cusps on
the two edges of each digon lies on different sides.   The definition
of $J$ guaranteed that the train track
$\bdd B$ on $T$ can be split to produce a curve $C$ of slope $\rrr$ on
$T$. Split $\bdd B \times I$ accordingly and, after Dehn filling,  cap
off $C$ by a meridian disk in the Dehn filling solid torus, one
obtains a branched surface $B_{\rrr}$ in $K(\rrr)$.  It was shown in
[Ro] that $B'$ carries an essential lamination in $K(\rrr)$.

Denote by $E(B)$ the exterior of $B$ in $E(K)$, i.e $E(B) = E(K) 
- \Int N(B) $.

\proclaim{Corollary 5.2}
For any $\rrr \in J$, the manifold $E(B)$ is homeomorphic to a 
component $W$ of the exterior of $B_{\rrr}$ in $K(\rrr)$, with
horizontal surface of $B$ identified to the horizontal surface of
$B_{\rrr}$ on $\bdd W$.
\endproclaim

\proof
Examine the above construction.  After adding the digons of
$T\times 1 - B$ to $B$, the branched surface is topologically
homeomorphic to $B \cup (T \times 1)$, which cuts off a region
isotopic to $E(B)$.  Clearly this region is not affected by
the later changes, and its horizontal surface is the
restriction of that on $E(B)$.
\qed
\enddemo

\head 6.  Surgery on twist knots, and the final proof
\endhead

As noticed earlier, any twist knot can be written as $K_{[2n, \pm
2]}$ for some $n$, because if $b$ is odd then $K_{[b, 2]} = K_{[b-1,
-2]}$.  Since $K_{[2n, -2]}$ is the mirror image of $K_{[-2n, 2]}$,
we need only consider knots of type $K_{[2n, 2]}$.  We can also
assume that $n \neq 0, 1$, otherwise the knot is a trivial knot or a
trefoil knot.

\proclaim{Lemma 6.1} Let $K = K_{[2n, 2]}$ be a twist knot with $n
\neq 0, 1$.  Then $K(\rrr)$ is not a small Seifert fiber space for
all $\rrr < -4$.  
\endproclaim

\proof
A knot $K = K_{[2n, 2]}$ has two spanning surfaces as indicated in
Figure 8, where $2n = 4$.  The first surface is a punctured Klein
bottle, and the second one is a punctured torus.

\bigskip

\bigskip
\leavevmode

\epsfxsize=5in
\centerline{\epsfbox{Fig8.ai}}
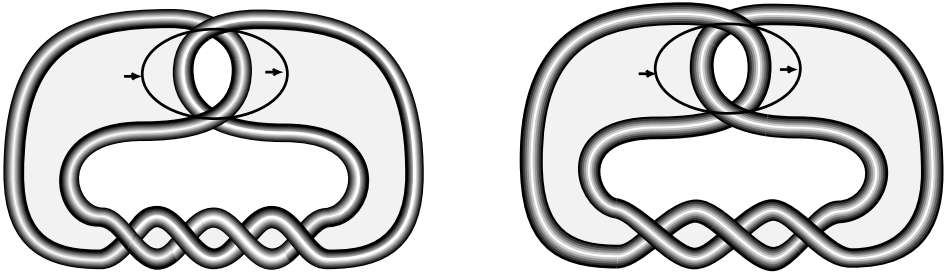
\centerline{Figure 8}

\bigskip

Let $R$ be the punctured Klein bottle of Figure 8(1) in the knot
exterior.  Add a disk $S$ to $R$ such that the boundary of $S$ is
the circle indicated in Figure 8(1).  The boundary of $S$ 
consists of four arcs $a_1 \cup b_1 \cup a_2 \cup b_2$, where
$b_i$ lies on $R$, and $a_i$ on the torus $T = \bdd N(K)$.  The
arrows at the arcs $b_i$ indicate the side of the cusp.  This
determines the branched surface $B = \left< R, S\right>$, and
hence the type of $a_i$ on $T$.  One of the $a_i$ (the one on top) is
of type $P_1$, and the other one of type $N_2$.  Using the notation in
Section 5, we have $p_1 = n_2 = 1$, and $p_2=n_1= 0$.  By calculating
the linking number between $K$ and $\bdd R$, one can see that $\bdd
R$ has slope $r = -4$ on $T$.  Therefore, 
 $$ \align
J  & = \{ r + (p_1 - n_1) \frac {x}{x+1} + (p_2 -
n_2) x \, \, \vert \, \, x > 0 \} = \{-4 +\frac{x}{x+1} - x \,\,
\vert \,\, x>0 \} \\
& = \{ -4 - \frac{x^2}{x+1} \,\, \vert \,\, x>0\}
= (-\infty, -4)
\endalign
$$

We need to show that $B$ is an essential branched surface in
$E(K)$.  It is clear that $N(B)$ is topologically a solid torus.  By
examining the cusp on $\bdd N(B)$, one can see that the exterior of
$B$ is the same as that of the surface in Figure 9, with cusp the
boundary of the surface $F$.  Note that $\bdd F$ runs along the
meridian of the solid torus $N(F)$ $2n-1$ times.  Hence $E(B)$, the
exterior of $B$,  is a solid torus with a single cusp running along
the longitude $2n-1$ times.  Since $n \neq 0, 1$, it is easy to see
that $E(B)$ satisfies all conditions for $B$ to be essential, i.e,
$E(B)$ is indecomposable, and the horizontal surface on $\bdd E(B)$ is
essential.  One also needs to check that the branched surface $B$
satisfies all the intrinsic essentiality properties, i.e, it has no
disk or half disk of contact, it contains no Reeb branched surfaces,
and it fully carries a lamination.  This is straight
forward.

\bigskip

\bigskip
\leavevmode

\epsfxsize=2.5in
\centerline{\epsfbox{Fig9.ai}}
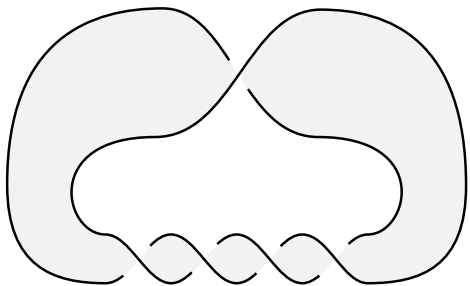
\centerline{Figure 9}

\bigskip

To apply Roberts' theorem, we also need to show that $B= \left<R,
S\right>$ has no planar surface of contact.  Assume that $P$ is a
compact surface embedded in $N(B)$, transverse to the $I$-fibers,
with $\bdd P \subset \bdd_vN(B) \cup T$.  Cutting along the two arcs
$b_1, b_2$ along which $S$ is glued to $R$, the branched surface $B$
becomes two surface $S$ and $R' = R - b_1\cup b_2$.  Let $u, v$ be
the number of times that $P$ intersects an $I$-fiber of $S$ and $R'$
respectively.  Then the gluing at $b_i$ gives the equation $v
= u+v + w$, where $w$ is the times $\bdd P$ passes the cusp at
$b_i$.  So $u=w =0$, and $P$ lies in $N(R)$.  Since $R$ is nonplanar,
$P$ can not be planar.  

It now follows from Theorem 5.1 that $B$ extends to an essential
branched surface $B_{\rrr}$ for all slopes $\rrr < -4$.  By
Corollary 5.2, $E(B)$ can be considered as a component of the
exterior of $B_{\rrr}$.  Since $E(B)$ is a solid torus with a cusp
running along the longitude $|2n-1| \geq 3$ times, it is not an
$I$-bundle.  Therefore, by the result of Brittenham [Br], $K(\rrr)$
is not a small Seifert fiber space for all $\rrr < -4$.  
\qed
\enddemo

\proclaim{Lemma 6.2} Let $K = K_{[2n, 2]}$ be a twist knot with $|n|
> 2$.  Then $K(\rrr)$ is not a small Seifert fiber space for
all $\rrr > 0$.  
\endproclaim

\proof
The proof is very similar to that of Lemma 6.1, only that instead of
using the non orientable surface in Figure 8(1), we use the orientable
surface $R$ in Figure 8(2).  As a Seifert surface of $K$, the boundary
of $R$ has slope $0$.  Also, the arcs $a_1, a_2$ are now of type
$P_2$ and $N_1$, so $p_1=n_2 = 0$, $p_2=n_1=1$, and
 $$ \align
J  & = \{ r + (p_1 - n_1) \frac {x}{x+1} + (p_2 -
n_2) x \, \, \vert \, \, x > 0 \} = \{0 -\frac{x}{x+1} + x \,\,
\vert \,\, x>0 \} \\
& = \{\frac{x^2}{x+1} \,\, \vert \,\, x>0\}
= (0, \infty)
\endalign
$$

The exterior of the branched surface $B = \left<R, S\right>$ is the
same as that of a band with $2n$ twists.  Since $|n| > 2$, it is not
an $I$-bundle.  Therefore, one can use the argument in the proof of
Lemma 6.1 to obtain a conclusion.
\qed
\enddemo

\bigskip

\demo{Proof of Theorem 1.1} Parts (1) and (2) are exactly Corollaries
3.4 and 4.2.  For part (3), consider the knot $K_{[2n, 2]}$.  By Lemma
2.2 and Corollary 2.4, surgeries with $ \rrr =0, 4$ are the only
toroidal ones, and $\rrr = -1, -2, -3$ are Seifert fibered.  By
Corollary 4.3, Lemmas 6.1 and 6.2, there are no other small Seifert
fibered surgeries.  Since no surgeries on $K_{[2n, 2]}$ are reducible
(Lemma 2.1), all surgeries with $\rrr \neq 0, -1, -2, -3, -4$ are
non-exceptional.  Since $K_{[2n, -2]}$ is the mirror image of
$K_{[-2n, 2]}$, the result follows.

Part (4) follows from the well known result of Thurston about surgery
on the Figure 8 knot [Th1], which is stronger as the non-exceptional
surgeries are shown to be hyperbolic.  The result as stated can also
be proved using the techniques in this paper: As toroidal and
reducible surgeries are known, we only need to deal with small Seifert
fibered surgeries.  By [Br] all non-integral surgeries on $K = K_{[-2,
2]}$ are not small Seifert fibered.  By Lemma 6.1 $K(\rrr)$ are not
small Seifert fibered for $\rrr < -4$.  Since $K$ is amphicheiral,
this is also true for $\rrr > 4$.  \qed \enddemo

\remark{Remark 6.3}  
Boileau and Porti [BP] proved a version of Thurston's Orbifold
Theorem, showing that the geometrization conjecture is true for
manifolds which admit a finite group action with nonempty fixed point
set.  Using this result one can show that non-exceptional surgeries on
2-bridge knots are also hyperbolic.  The proof is as follows.  A $p/q$
2-bridge knot can be obtained by taking two arcs of slope $p/q$ on the
``pillowcase'', then joining the ends with two trivial arcs.  From
this picture it is easy to see that $K$ is a strongly invertible knot,
i.e.\ there is an involution $\varphi$ of $S^3$ such that $\varphi(K)
= K$, and the fixed point set of $\varphi$ is a circle $S$
intersecting $K$ at two points.  The map $\varphi$ restricts to an
involution of $E(K) = S^3 - \Int N(K)$, which can be extended to an
involution $\hat{\varphi}$ of the surgered manifold.  Since $\hat
{\varphi}$ has nonempty fixed point set, the result follows from [BP].
\endremark

\enddocument